\documentclass[11pt]{amsart}
\usepackage{latexsym,amssymb,amsthm,amsmath,amsfonts,amsxtra,graphicx}

\def\R{\mathbb{R}}
\def\H{\mathbb{H}}

\def\dt{\partial_t}
\def\du{\partial_u}
\def\dv{\partial_v}

\def\E{\mathbb{E}}
\def\S{\mathbb{S}}
\def\({\left(}
\def\){\right)}

\theoremstyle{plain}

\newtheorem{theorem}{Theorem}

\newtheorem{proposition}{Proposition}

\theoremstyle{remark}
\newtheorem{remark}{Remark}

\title{Constant Angle Surfaces in a warped product}

\author[F. Dillen]{Franki Dillen}
\address[F. Dillen]{Katholieke Universiteit Leuven\\ Departement
Wiskunde\\ Celestij\-nenlaan 200 B, Box 2400\\ BE-3001 Leuven\\ Belgium}
\email[F. Dillen]{franki.dillen@wis.kuleuven.be}

\author[M.I. Munteanu]{Marian Ioan Munteanu}
\address[M.I. Munteanu]{University 'Al.I.Cuza' of Ia\c si\\ Faculty of Mathematics\\
Bd. Carol I, no.11\\ 700506 Ia\c si\\ Romania} \email[M.I.
Munteanu]{marian.ioan.munteanu@gmail.com}

\author[J. Van der Veken]{Joeri Van der Veken}
\address[J. Van der Veken]{Katholieke Universiteit Leuven\\ Departement
Wiskunde\\ Celes\-tijnen\-laan 200 B, Box 2400\\ BE-3001 Leuven\\ Belgium}
\email[J. Van der Veken]{joeri.vanderveken@wis.kuleuven.be}

\author[L. Vrancken]{Luc Vrancken}
\address[L. Vrancken]{Univ. Lille Nord de France\\ F-59000 Lille\\
France; UVHC\\ LAMAV\\ F-59313 Valenciennes\\ France; Katholieke
Universiteit Leuven\\ Departement Wiskunde\\ Celes\-tijnen\-laan 200 B, Box 2400\\ BE-3001 Leuven\\
Belgium}\email[L. Vrancken]{luc.vrancken@univ-valenciennes.fr}

\date{\today }

\thanks{The second author was supported by Grant PN-II ID 398/2007-2010 (Romania)}
\thanks{The third author is a postdoctoral researcher of the Research Foundation - Flanders (F.W.O.)}
\thanks{Research supported also by Research Foundation - Flanders project G.0432.07.}

\begin{document}

%\input epsfx.tex

%\pagestyle{plain}
%\pagenumbering{arabic}

\begin{abstract}
Let $I \subseteq \R$ be an open interval, $f : I \to \R$ a
strictly positive function and denote by $\E^2$ the Euclidean plane. We classify all surfaces in the warped
product manifold $I \times_f \E^2$ for which the unit normal makes a
constant angle with the direction tangent to $I$.
\end{abstract}



\maketitle

\section{Introduction}

In the last few years, the study of the geometry of surfaces in
3-dimensional spaces, in particular of product type $M^2\times\R$
was developed by a large number of mathematicians.
In particular, in \cite{kn:DFVV07}, \cite{kn:DM07} and
\cite{kn:DM09} the authors have studied constant angle surfaces in
$S^2\times \R$ and $H^2\times \R$, namely those surfaces for which
the unit normal makes a constant angle with the tangent direction to
$\R$. In \cite{kn:FMV08} a classification of surfaces in the
3-dimensional Heisenberg group making a constant angle with the
fibers of the Hopf-fibration was obtained. In all these spaces, the
angle which is required to be constant is one of the fundamental
invariants appearing in the existence and uniqueness theorem for
isometric immersions, cfr. \cite{daniel}. In another recent paper
\cite{kn:CdS07} it is proven that if the ambient space is the Euclidean
3-space, the study of surfaces making a constant angle with a fixed
direction has some important applications to physics, namely in
special equilibrium configurations of nematic and smectic C liquid
crystals. In \cite{kn:LM09} constant angle surfaces in
3-dimensional Minkowski space were studied.

In the present paper, constant angle surfaces in another important
family of 3-spaces in which there exists a distinct direction,
namely warped products of an open interval with the Euclidean plane,
are classified. Special examples, such as flat or minimal surfaces
in this family are given.

\section{Preliminaries}
The following notion of warped product or, more generally, warped
bundle was introduced by Bishop and O'Neill in \cite{kn:BN69} in
order to construct a large variety of manifolds of negative
curvature. Let $B$ and $F$ be two Riemannian manifolds with
Riemannian metrics $g_B$ and $g_F$ respectively. Let $f>0$ be a
smooth positive function on $B$ and denote by $B\times F$ the
product manifold. The warped product of $B$ and $F$ with warping
function $f$ is the Riemannian manifold
\begin{equation*}
B \times_f F = \(B\times F,\ g_B+f^2\,g_F\).
\end{equation*}

Let $f:I\subseteq\R\to\R$ be a smooth strictly positive function on an open interval $I$ and consider the
warped product of $I$ and the Euclidean plane $\E^2$ with warping function $f$
\begin{equation*}
\label{2.2} (\widetilde M, \widetilde g) = I \times_f \E^2 = \(I
\times \R^2,\ dt^2+f(t)^2(dx^2+dy^2)\)
\end{equation*}
where $t$ is the coordinate on $I$ and $x$ and $y$ are coordinates
on $\E^2$.

Denote by $\widetilde \nabla$ the Levi-Civita connection of
$(\widetilde M, \widetilde g)$. Denote by $U$, $V$ and $W$ lifts of
vector fields tangent to $\E^2$. One has
\begin{subequations}
\renewcommand{\theequation}{\theparentequation .\alph{equation}}
\label{eq:1}
\begin{align}
\label{2.3} & \widetilde \nabla_U V = D_U V - \frac{f'}{f}~\widetilde g (U,V)~\dt\\
\label{2.4} & \widetilde \nabla_U \dt = \widetilde \nabla_{\dt} U = \frac{f'}{f}~U\\
\label{2.5} & \widetilde \nabla_{\dt} \dt = 0
\end{align}
\end{subequations}
where $D$ is the covariant derivative on $\E^2$, see for
example \cite{ON}. Remark that we have identified $U$ and $V$ with
their projections onto $\E^2$. From these equations, it follows
immediately that the curvature tensor $\widetilde R$, defined as
$\widetilde R(U,V)=[\widetilde\nabla_U,\widetilde\nabla_V]-\widetilde\nabla_{[U,V]}$
is given by
\begin{subequations}
\renewcommand{\theequation}{\theparentequation .\alph{equation}}
\label{eq:2}
\begin{align}
\label{2.6} & \widetilde R(U,\dt)V = \frac{f''}{f}~\widetilde g(U,V)~\dt\\
\label{2.7} & \widetilde R(U,V)\dt = 0\\
\label{2.8} & \widetilde R(U,\dt)\dt = -\frac{f''}{f}~U\\
\label{2.9} & \widetilde R(U,V)W = -\frac{(f')^2}{f^2}~\big(\widetilde g(V,W) U -\widetilde g(U,W) V\big).
\end{align}
\end{subequations}
Let $\iota : M \to \widetilde M$ be an immersion of a surface $M$ in
$\widetilde M$ and let $g$ be the pulled back metric of $\widetilde
g$ on $M$. We will not write down $\iota$, unless it is absolutely necessary to avoid confusion.
Let $\xi$ be a unit normal vector field on $M$ and denote
by $A$ the associated shape operator. The formulas of Gauss and
Weingarten state respectively that\\
$({\mathbf{G}})
%\label{2.10}
\qquad\qquad
\widetilde\nabla_X Y = \nabla_X Y + h(X,Y)$\\
$({\mathbf{W}})
%\label{2.11}
\qquad\qquad
\widetilde\nabla_X \xi = -AX$\\
for every $X$ and $Y$ tangent to $M$. Here, $\nabla$ is the
Levi-Civita connection of $M$ and $h$ is the second
fundamental form. We have $\widetilde g(h(X,Y),\xi) = g(X,AY)$ for
all $X$ and $Y$ tangent to $M$. One can decompose $\dt$ as
\begin{equation}
\label{2.12} \dt = T + \cos\theta\,\xi,
\end{equation}
where $\theta\in[0,\pi)$ is the angle between $\dt$ and the normal $\xi$ and
$T$ is the projection of $\dt$ on the tangent plane of $M$. We have
$\cos\theta=\widetilde g(\dt,\xi)$ and, since $\dt$ has unit length,
$|T|=\sin\theta$.

If one denotes by $R$ the curvature tensor on $M$, then it follows
from \eqref{eq:2}, (G), (W) and \eqref{2.12} that the equations of Gauss and
Codazzi can be written respectively as\\[2mm]
$({\mathbf{EG}})
\begin{array}{rl}
%\label{2.13}
\quad
R(X,Y)Z = & \ g(AY,Z)AX - g(AX,Z)AY \\
& - \((\log f)'\circ\iota\)^2 (g(Y,Z)X - g(X,Z)Y) \\
& - \((\log f)''\circ\iota\) \Big(g(Y,T)g(Z,T)X - g(X,T)g(Z,T)Y  \\
& \hspace{2.7cm}  - g(Y,T)g(X,Z)T + g(X,T)g(Y,Z)T \Big)
\end{array}
$\\[2mm]
$
({\mathbf{EC}})
\qquad
%\label{2.14}
(\nabla_X A) Y-(\nabla_Y A)X = \cos\theta \( (\log
f)''\circ\iota \) \(g(Y,T)X - g(X,T)Y\)
$\\
for $X$, $Y$ and $Z$ tangent to $M$.
\begin{proposition}
\label{prop2.1}
Let $X$ be tangent to $M$, then
\begin{align}
\label{2.15} & \nabla _XT = \cos\theta\,AX + \((\log f)'\circ \iota\)\(X-g(X,T)T\),\\
\label{2.16} & X(\cos\theta) = -g(X,AT)-\cos\theta \((\log f)'\circ\iota\) g(X,T).
\end{align}
\end{proposition}
\proof If $X$ is tangent to $M$, then $\widetilde g(X,\dt)=g(X,T)$.
One can express $\widetilde \nabla _X\dt$ in two ways:

\quad $\widetilde \nabla _X\dt = \((\log f)'\circ\iota\) \left(
X-g(X,T)\dt\right)$, by use of \eqref{2.4} and \eqref{2.5},

\quad $\widetilde\nabla_X \dt = \nabla_X T + h(X,T) +
X(\cos\theta)\xi - \cos\theta AX$, by use of (G), (W) and \eqref{2.12}. \\
Comparing the tangent and the normal parts, one gets the conclusion.
\endproof

From \eqref{2.16} we obtain immediately the following.
\begin{proposition}
\label{prop2.2}
If $\theta$ is a constant angle, then $T$ is a principal direction
and the corresponding eigenvalue of the shape operator is
$-\cos\theta\((\log f)'\circ\iota\)$.
\end{proposition}

From now on, we will assume that $\theta$ is constant. In this case we say that
$\iota:M\to\widetilde M$ is {\em a constant angle
surface}.

We may assume that $\theta\in[0,\pi/2]$.

If $\theta=0$, then $\iota(M)\subseteq  \{t_0\}\times \E^2$, so we suppose that $\theta\neq 0$.
Then $T\neq 0$ and we can consider $e_1=T/|T|=T/\sin\theta$. Let
$e_2$ be a unit tangent vector, orthogonal to $e_1$. Then $e_2$ is also a principal direction, thus there exists a function
$\lambda\in C^\infty(M)$ such that $Ae_2=\lambda e_2$. Combining with \eqref{2.15}, this yields the following.
\begin{proposition} \label{prop2.3}
Let $M$ be a constant angle surface in $\widetilde M$, with
$\theta\neq 0$. Then there exists an orthonormal frame field
$\{e_1,e_2\}$ on $M$ such that the shape operator with respect to
this frame takes the form
\begin{equation}
\label{2.17}
A=\left(\begin{array}{cc}
-\cos\theta\((\log f)'\circ\iota\) & 0\\[2mm]
0 & \lambda
\end{array}\right)
\end{equation}
for some $\lambda\in C^\infty(M)$ and the Levi-Civita connection is
given by
\begin{equation}
\begin{aligned}
\label{2.18} & \nabla_{e_1}e_1 = 0, & & \qquad \nabla_{e_2}e_1 =
        \frac{1}{\sin\theta}\(\lambda\cos\theta+\((\log f)'\circ\iota\)\)e_2, \\
& \nabla_{e_1}e_2 = 0,  & & \qquad \nabla_{e_2}e_2 =
-\frac{1}{{\sin\theta}}\(\lambda\cos\theta+\((\log f)'\circ\iota\)\)
e_1.
\end{aligned}
\end{equation}
\end{proposition}

\section{The classification theorem}

In this section we classify the constant angle surfaces in
$(\widetilde M, \widetilde g) = I \times_f \E^2$ with
$\theta\neq 0$. We consider the orthonormal frame field $\{e_1,e_2\}$ as
above. Then from \eqref{2.18} we obtain that
$[e_1,e_2]$ is proportional to $e_2$. Therefore we can
choose coordinates $(u,v)$ such that $\partial_u=e_1$ and
$\partial_v=\beta e_2$ for some function $\beta$. Then it is clear that $g$ takes the form
\begin{equation}
\label{3.2} g=du^2+\beta^2(u,v)\ dv^2.
\end{equation}
The Levi-Civita connection is determined by
\begin{equation}
\label{3.3} \nabla_{\du}\du=0, \ \
\nabla_{\du}\dv=\nabla_{\dv}\du=\frac{\beta_u}\beta\,\dv,\ \
\nabla_{\dv}\dv=-\beta\beta_u\du+\frac{\beta_v}\beta\,\dv
\end{equation}
and $\beta$ satisfies
\begin{equation}
\label{3.1}
\beta_u = \frac{\beta}{\sin\theta} \(
\lambda\cos\theta +\((\log f)'\circ\iota\)\).
\end{equation}
If we put
$$
\iota(u,v)=(t(u,v),x(u,v),y(u,v))
$$
then
$$
t_u = \widetilde g (\iota_u,\partial_t) = \widetilde g
(e_1,\partial_t) = \widetilde g (T/\sin\theta,T+\cos\theta\xi) =
\sin\theta
$$
and
$$
t_v = \widetilde g (\iota_v,\partial_t) = \widetilde g
(\beta e_2,\partial_t) = \widetilde g (\beta e_2,T+\cos\theta\xi) = 0
$$
such that, after a translation in the $u$ coordinate,
\begin{equation}
\label{3.4} t(u,v) = u\sin\theta.
\end{equation}

\begin{theorem} \label{theo3.1}
An isometric immersion $\iota:M\to I\times_f\E^2$ defines a surface with constant
angle $\theta\in [0,\pi/2]$ if and only if, up to rigid
motions of $I\times_f\E^2$, one of the following holds locally.
\medskip

\noindent $\mathrm{(i)}$ There exist local coordinates $(u,v)$ on
$M$, with respect to which the immersion $\iota$ is given by
\begin{multline}
\label{3.13} \iota(u,v) = \( u\sin\theta,\
\cot\theta \( \int^{u\sin\theta} \frac{d\tau}{f(\tau)} \) \cos v - \int^v \alpha(\tau)\sin\tau d\tau, \right.\\
\left. \cot\theta \( \int^{u\sin\theta} \frac{d\tau}{f(\tau)} \)
\sin v + \int^v \alpha(\tau)\cos\tau d\tau \)
\end{multline}
for some smooth function $\alpha$.
\medskip

\noindent $\mathrm{(ii)}$ $\iota(M)$ is an open part of the cylinder
\begin{equation} \label{3.14}
x - G(t) = 0,
\end{equation}
where $G(t)=\cot\theta \int^t\frac{d\tau}{f(\tau)}$. This surface is totally umbilical with mean curvature $H=-\cos\theta f'(u\sin\theta)/f(u\sin\theta)$.

\noindent $\mathrm{(iii)}$ $\iota(M)$ is an open part of the surface
$t=t_0$ for some real number $t_0$, and $\theta=0$.
\end{theorem}

\proof
Let us first check that the surfaces described in the theorem are constant angle surfaces.

For case (i), a basis for the tangent plane to the surface is given
by
\begin{align*}
 & \iota_u = \( \sin\theta,\ \frac{\cos\theta\cos v}{f(u\sin\theta)},\ \frac{\cos\theta\sin v}{f(u\sin\theta)} \right)\\
 & \iota_v = \( \cot\theta \( \int^{u\sin\theta}
\frac{d\tau}{f(\tau)} \) +\alpha(v) \) \( 0, -\sin v,  \cos v \).
\end{align*}
Notice that if $a=(a_1,a_2,a_3)$ and $b=(b_1,b_2,b_3)$ are vectors
in $T_{(t,x,y)}(I\times_f\E^2)$, then the vector defined by
\begin{equation*}
 a \times_f b = \( f^2(t)(a_2 b_3 - a_3 b_2),\ a_3 b_1 -
a_1 b_3,\ a_1 b_2 - a_2 b_1 \)
\end{equation*}
is orthogonal to both $a$ and $b$. Hence
\begin{equation*}
\xi = \frac{\iota_u \times_f \iota_v}{|\iota_u \times_f
\iota_v|} = \( \cos\theta,\ -\frac{\sin\theta\cos
v}{f(u\sin\theta)},\ -\frac{\sin\theta\sin v}{f(u\sin\theta)}\right)
\end{equation*}
is a unit normal on the surface. We immediately deduce that $\widetilde
g(\xi,\dt) = \cos\theta$.

For case (ii), one can use the parametrization
\begin{equation*}
\iota(u,v) = \( u,\ \cot\theta \int^u \frac{d\tau}{f(\tau)},\ v \).
\end{equation*}
Then $\xi=(\cos\theta, -\sin\theta/f(u),0)$ is a unit normal and $\widetilde
g(\xi,\dt) = \cos\theta$.

Case (iii) is obvious.
\medskip

Conversely, let $\iota:M\to I\times_f\E^2$ be a constant angle
surface with constant angle $\theta$. As mentioned before, we may
assume that $\theta\in[0,\pi/2]$. If $\theta=0$ then $\iota(M)$ is of type (iii) described in the theorem.
If $\theta=\pi/2$, the vector field $\partial_t$ is everywhere
tangent to $\iota(M)$. This implies that $\iota(M)$ is an open part
of a cylinder with rulings in the direction of $\partial_t$ or,
equivalently that there exist local coordinates $(u,v)$ on $M$ such
that $\iota(u,v)=(u,\gamma_1(v),\gamma_2(v))$ for some smooth
functions $\gamma_1$ and $\gamma_2$. If $\iota$ parametrizes a
plane, this is case (ii) of the theorem with $\theta=\pi/2$. If
$\iota$ does not describe a plane, this is case (i) of the theorem
with $\theta=\pi/2$.

From now on, assume that $\theta\in (0,\pi/2)$. If we choose local
coordinates on $M$ as above, we can write $\iota(u,v)=(u\sin\theta,x(u,v),y(u,v))$.
Using \eqref{3.2} we obtain
\begin{subequations}
\renewcommand{\theequation}{\theparentequation .\alph{equation}}
\label{eq:23}
\begin{align}
\label{3.22} & f^2(u\sin\theta)\(x_u^2+y_u^2\) = \cos^2\theta \\
\label{3.23} & x_ux_v + y_uy_v = 0 \\
\label{3.24} & f^2(u\sin\theta)\(x_v^2+y_v^2\) = \beta^2.
\end{align}
\end{subequations}
Define
\begin{equation}
\label{3.25} \sigma(u) = \log f(u\sin\theta) = ((\log
f)\circ\iota)(u,v).
\end{equation}
Then a straightforward computation, using
\eqref{eq:1}, \eqref{eq:23} and \eqref{3.25} yields
\begin{subequations}
\renewcommand{\theequation}{\theparentequation .\alph{equation}}
\label{eq:27}
\begin{align}
\label{3.26} & \widetilde \nabla_{\iota_u}\iota_u
= \iota_{uu} + 2\sigma' \iota_u - \(\sin\theta+\frac1{\sin\theta}\)\sigma'\dt \\
\label{3.27} & \widetilde \nabla_{\iota_u}\iota_v = \iota_{uv} + \sigma' \iota_v \\
\label{3.28} & \widetilde\nabla_{\iota_v}\iota_v = \iota_{vv} -
\frac{1}{\sin\theta}\beta^2\sigma'\dt.
\end{align}
\end{subequations}
On the other hand, we can express these covariant derivatives by
using the formula of Gauss (G). By using \eqref{2.12},
\eqref{2.17}, \eqref{3.3} and \eqref{3.25} we obtain
\begin{subequations}
\renewcommand{\theequation}{\theparentequation .\alph{equation}}
\label{eq:30}
\begin{align}
\label{3.29} & \widetilde \nabla_{\iota_u}\iota_u = \sigma' \iota_u - \frac{1}{\sin\theta}\sigma'\dt, \\
\label{3.30} & \widetilde \nabla_{\iota_u}\iota_v = \frac{\beta_u}{\beta} \iota_v, \\
\label{3.31} & \widetilde\nabla_{\iota_v}\iota_v =
-\(\beta\beta_u+\tan\theta\lambda\beta^2\)\iota_u +
\frac{\beta_v}{\beta}\iota_v +
\frac{1}{\cos\theta}\lambda\beta^2\partial_t.
\end{align}
\end{subequations}
We will now compare successively \eqref{eq:27} to \eqref{eq:30}.\\
From \eqref{3.26} and \eqref{3.29} we obtain
\begin{equation*}
\iota_{uu} + \sigma' \iota_u -
\sin\theta\sigma'\partial_t = 0.
\end{equation*}
This equation is satisfied for the $t$-component. For the $x$- and
the $y$-component we obtain respectively $x_{uu}+\sigma' x_u=0$ and
$y_{uu}+\sigma' y_u=0$, such that $x_u(u,v)=e^{-\sigma(u)}c_1(v)$
and $y_u(u,v)=e^{-\sigma(u)}c_2(v)$ for some functions $c_1$ and
$c_2$. From \eqref{3.22} we obtain $c_1^2(v)+c_2^2(v)=\cos^2\theta$.
If we put $p_1(v)=c_1(v)/\cos\theta$ and $p_2(v)=c_2(v)/\cos \theta$,
then
\begin{equation}
\label{3.33} \iota_u(u,v) = \( \sin\theta,\ \cos\theta
e^{-\sigma(u)} p_1(v),\ \cos\theta e^{-\sigma(u)} p_2(v) \), \ \
p_1^2(v)+p_2^2(v)=1.
\end{equation}
From \eqref{3.27} and \eqref{3.30}, we obtain
\begin{equation*}
\iota_{uv} + \(\sigma'-\frac{\beta_u}{\beta}\) \iota_v
= 0.
\end{equation*}
This equation is again satisfied for the $t$-component. Integrating, we obtain
\begin{equation}
\label{3.35} \iota_v(u,v) = e^{-\sigma(u)}\beta(u,v) \( 0,\ q_1(v),\
q_2(v) \), \ \ q_1^2(v)+q_2^2(v)=1.
\end{equation}
Remark that the compatibility condition for \eqref{3.33} and
\eqref{3.35} is
\begin{equation}
\label{3.36} (p_1',p_2') =
\frac{1}{\cos\theta}(\beta_u-\sigma'\beta)\(q_1,q_2\).
\end{equation}
Finally, from \eqref{3.28} and \eqref{3.31}, we obtain
\begin{equation}
\label{3.37} \iota_{vv} +
\(\beta\beta_u+\tan\theta\lambda\beta^2\)\iota_u  -
\frac{\beta_v}{\beta}\iota_v -
\beta^2\(\frac{\sigma'}{\sin\theta}+\frac{\lambda}{\cos\theta}\)\partial_t
= 0.
\end{equation}
If we substitute \eqref{3.33} and \eqref{3.35} into \eqref{3.37},
the resulting equations for the $x$- and the $y$-component yield
\begin{equation}
\label{3.38} (q_1',q_2') =
-(\beta_u\cos\theta+\lambda\beta\sin\theta) \(p_1,p_2\).
\end{equation}

At this point we can distinct two cases: $(p_1(v), p_2(v))$ is constant or not.

\medskip

\noindent\underline{Case 1}: $(p_1(v), p_2(v))$ is constant.

Then from \eqref{3.36}  we obtain that $\beta_u=\sigma' \beta$, and hence $\beta(u,v) = \psi(v) f(u\sin\theta)$.
After a change in the $v$-coordinate, we can assume that $\psi(v)=1$, such that $\beta(u,v)=f(u\sin\theta)$.
From \eqref{3.1} we then obtain that $\lambda = -\cos\theta f'(u\sin\theta)/f(u\sin\theta)$. From Proposition \ref{prop2.3}
it follows that $M$ is totally umbilical.

From \eqref{3.38} then follows that $(q_1,q_2)$ is constant. Integrating \eqref{3.33} and using \eqref{3.35}
gives us
\begin{multline}
\label{3.48} \iota(u,v) = \( u\sin\theta,\ p_1\cos\theta \(\int^u e^{-\sigma(\mu)}d\mu\) + q_1 v +a_1, \right. \\
\left. p_2\cos\theta \(\int^u e^{-\sigma(\mu)}d\mu\) +
q_2 v +a_2 \)
\end{multline}
for some constants $a_1$ and $a_2$, which can be taken zero after a translation in $x$ and $y$. Moreover, since $\widetilde g(\iota_u,\iota_v)=0$, we have $p_1q_1+p_2q_2=0$. Hence, after a rotation around the $t$-axis, which is an isometry of $I\times_f\E^2$, we may assume that $(p_1,p_2)=(1,0)$ and $(q_1,q_2)=(0,1)$. Hence we obtain after a substitution
$\tau=\mu\sin\theta$
\begin{equation*}
\iota(u,v) = \( u\sin\theta,\
\cot\theta\int^{u\sin\theta}\frac{d\tau}{f(\tau)},\ v \)
\end{equation*}
which corresponds to case (ii) of the theorem.

\smallskip
%%%%%%%%%%%%%%% CASE 2 %%%%%%%%%%%%%%%%%%%%

\noindent\underline{Case 2}: $(p_1(v), p_2(v))$ is not constant.
Then from \eqref{3.33} we can assume that, after a change of the $v$-coordinate, that
\begin{equation}\label{3.33b}
(p_1(v), p_2(v)) = (\cos v, \sin v).
\end{equation}
Then \eqref{3.36} implies that
\begin{equation}\label{3.33c}
\beta_u-\sigma'\beta = \pm \cos\theta
\end{equation}
and by changing the sign of $u$, we can assume the right hand side to be $\cos\theta$.
Integrating \eqref{3.33c} gives
\begin{equation}\label{3.33d}
\beta(u,v) e^{-\sigma(u)} - \cos\theta \int^u e^{-\sigma(\mu)}d\mu = \alpha(v)
\end{equation}
for some function $\alpha(v)$.
Furthermore \eqref{3.36}
shows that
\begin{equation*}
(q_1(v), q_2(v)) = (-\sin v,\cos v).
\end{equation*}
Hence \eqref{3.33} and \eqref{3.35} reduce to
\begin{align}
\label{3.42} & \iota_u(u,v) = \( \sin\theta,\ \cos\theta e^{-\sigma(u)}\cos v,\cos\theta e^{-\sigma(u)}\sin v \) \\
\label{3.43} & \iota_v(u,v) = e^{-\sigma(u)}\beta(u,v) \( 0, -\sin
v,\cos v \).
\end{align}
Integrating \eqref{3.42} gives
\begin{multline}
\label{3.44} \iota(u,v) = \( u\sin\theta,\ \cos\theta \(\int^u e^{-\sigma(\mu)}d\mu\) \cos v + \gamma_1(v),\right.\\
\left. \cos\theta \(\int^u e^{-\sigma(\mu)}d\mu\) \sin
v + \gamma_2(v) \)
\end{multline}
for some smooth functions $\gamma_1$ and $\gamma_2$.
If we take the derivative with respect to $v$ in
\eqref{3.44} and compare it to \eqref{3.43} we get, using \eqref{3.33d}
\begin{equation*}
( \gamma_1'(v),\gamma_2'(v) ) = \alpha(v)(-\sin v, \cos v).
\end{equation*}
After integration, we obtain case (i) of the theorem.
\begin{remark}
In this case, the function $\lambda$ is given by
\begin{equation}\label{lc2}
    \lambda\beta = \sin\theta - \frac{f'}{f}\beta \cos\theta.
\end{equation}
This follows from \eqref{3.1} and \eqref{3.38}.
\end{remark}
\begin{remark}
Notice that if we take the Euclidean metric on ${\mathbf{R}}^3$, i.e. the warping function is $1$,
we retrieve the statements of Theorem $7$ in \cite{kn:MN09}.
\end{remark}

\section{Rotational surfaces of constant angle}

In this section, we will classify constant angle surfaces in
$I\times_f\E^2$, which are invariant under rotations with respect
to the $t$-axis.

Let us first remark that any rotation
\begin{equation*}
R_{\phi} : I\times_f\E^2 \to I\times_f\E^2 : (t,x,y)
\mapsto (t,\ x\cos\phi-y\sin\phi,\ x\sin\phi+y\cos\phi)
\end{equation*}
is an isometry. Let $\gamma$ be a curve in the plane containing the
$t$- and the $x$-axis. Assume that $\gamma(u)=(a(u),b(u),0)$ is an
arc length parametrization, i.e., that
\begin{equation}\label{4.2}
(a'(u))^2 + f^2(a(u))(b'(u))^2 = 1.
\end{equation}
We want to investigate, under which conditions, the surface
\begin{equation*}
\iota(u,v) = \( a(u),\ b(u)\cos v,\ b(u)\sin v \)
\end{equation*}
is a constant angle surface in $I\times_f\E^2$.

The unit normal vector field is given by
\begin{equation*}
\xi(u,v) = \( b'(u)f(a(u)),\ -\frac{a'(u)\cos
v}{f(a(u))},\ -\frac{a'(u)\sin v}{f(a(u))} \).
\end{equation*}
Hence, the surface determines a constant angle surface
with constant angle $\theta$ if and only if
\begin{equation}
\label{4.5} b'(u)f(a(u))=\cos\theta.
\end{equation}
Combining \eqref{4.2} and \eqref{4.5} yields
\begin{equation}
\label{4.6} (a'(u))^2=\sin^2\theta.
\end{equation}
There are now two cases to consider.
\medskip

The case $\sin\theta=0$ is obvious and it corresponds to case (iii) of the Theorem \ref{theo3.1}.
So assume  $\sin\theta\neq 0$.
Then we see from \eqref{4.6} that $a(u)=\pm u\sin\theta+c$
for some real constant $c$. After a change of the arc length
parameter $u$ of $\gamma$, we may consider that
\begin{equation}
\label{4.7} a(u)=u\sin\theta.
\end{equation}
If $\theta=\pi/2$, then $b=b_0$ is constant and we obtain the circular cylinder
$\iota(u,v)=(u,b_0\cos v, b_0\sin v)$. In the sequel we will take $\theta\in (0,\pi/2)$.

It then follows from \eqref{4.5} that
\begin{equation*}
\label{4.8}
b(u) = \int^u \frac{\cos\theta}{f(\mu\sin\theta)}d\mu =
\cot\theta\int^{u\sin\theta}\frac{d\tau}{f(\tau)}.
\end{equation*}
We conclude that the rotational surface immersion becomes
\begin{multline}
\label{4.9} \iota(u,v) = \( u\sin\theta,\
\(\cot\theta\int^{u\sin\theta}\frac{d\tau}{f(\tau)}\)\cos v, \right. \\
\left.\(\cot\theta\int^{u\sin\theta}\frac{d\tau}{f(\tau)}\)\sin v \)
\end{multline}
which corresponds, up to a translation in the $x$-direction, to a
special case of case (i) of Theorem \ref{theo3.1}, namely the case
where $\alpha(v)=0$.

\section{Examples}

\subsection{Flat constant angle surfaces.}
\

\noindent A surface of type (iii) of Theorem \ref{theo3.1} is a trivial example of a flat surface with constant angle $\theta=0$.
In order to give an example of flat constant angle surface with
$\theta \neq 0$ consider a surface of type (ii) in Theorem
\ref{theo3.1}. Using (EG) and \eqref{2.17}, we obtain
\begin{equation*}
\label{eq:K} K = \det A - \left( (\log f)' \circ \iota \right)^2 -
\left( (\log f)'' \circ \iota \right) \sin^2\theta = - \left(
\frac{f''}{f} \circ \iota \right) \sin^2 \theta.
\end{equation*}
Thus $f(t)=a(t+b)$, with $a\ne 0$. The metric $\widetilde g$ on the
ambient space is called a \emph{cone metric}.

\subsection{Minimal constant angle surfaces.}
\

\noindent Consider first a constant angle surface of type (iii) of
Theorem \ref{theo3.1}. Then $\partial_t$ is a unit normal and it
follows from \eqref{2.4} that the surface is totally umbilical with
shape operator $A = f'(t_0)/f(t_0)\, \mathrm{id}$. Hence, such a
surface is minimal if and only if $f'(t_0)=0$, case in which it is
also totally geodesic.

Now assume that the constant angle surface is of  type (ii) of Theorem \ref{theo3.1}.
Then it is minimal only if it is totally geodesic. Since $
H = - \cos\theta {f'(u \sin\theta)}/{f(u \sin\theta)}$, either $\theta = \pi/2$, i.e. the surface is a warped
product of an interval and a straight line in $\E^2$, or $f'=0$,
i.e. the ambient space is a direct product and $M$ is a plane.

Finally, if we assume that the constant angle surface is of  type (i) of Theorem \ref{theo3.1},
then from \eqref{lc2} follows that $H=0$ if and only if
\begin{equation}\label{min}
 2\cos\theta \beta \sigma' = \sin^2\theta.
\end{equation}
Hence $\beta$ depends only on $u$. Differentiating \eqref{min} using \eqref{3.33c} yields
\begin{equation}
\label{min2}
\(\frac{1}{\sigma'}\)' = \frac{1+\cos^2\theta}{\sin^2\theta}.
\end{equation}
Integrating \eqref{min2} shows that $f$ has to take the form
\begin{equation*}
\label{min3}
    f(t)=b(t+c)^{\frac{\sin^2\theta}{1+\cos^2\theta}}.
\end{equation*}
Without loss of generality, we can assume $b=1$ and $c=0$. We put $m=\frac{\sin^2\theta}{1+\cos^2\theta}$, such that $f(t)=t^m$, $m\in(0,1)$.
From \eqref{prop2.3} and  \eqref{lc2} we then obtain that
$\lambda=\frac{m\cot\theta}{u}$ and $\beta=\frac{\cos\theta}{1-m}~u$. Then it follows that in \eqref{3.33d} we have to take $\alpha=0$.
Then from the classification Theorem \ref{theo3.1} we obtain that
\begin{equation*}
\iota(u,v) = \left( u \sin\theta, \frac{\cot\theta}{1-m} (u\sin \theta)^{1-m} \cos v, \frac{\cot\theta}{1-m} (u\sin \theta)^{1-m}\sin v \right).
\end{equation*}
This represents a constant angle minimal surface, with
$\theta=\arccos\sqrt{(1-m)/(1+m)}$.  Moreover, the
surface is a rotation surface.

\subsection{Constant angle surfaces with a harmonic height function.}
\

\noindent Consider the \emph{height function}
\begin{equation*}
\widetilde h: I\times_f\E^2 \to \R: (t,x,y) \mapsto t.
\end{equation*}
If $\iota: M\to I\times_f\E^2$ is an isometric immersion of a surface, then we denote by $h$
the restriction of $\widetilde h$ to $M$, i.e. $h = \widetilde g
(\iota, \dt)$. Remark that
$$
g(X,\mathrm{grad}\,h) = X(h) = X (\widetilde g(\iota,\dt)) = \widetilde g(X,\dt) = g(X,T)
$$
for all $X$ tangent to $M$ and hence
\begin{equation*}
\mathrm{grad}\,h = T.
\end{equation*}
Thus, by using \eqref{2.15} we obtain
\begin{equation}
\label{eq:laplacian}
\Delta h = \mathrm{div}\,T = \mathrm{trace}(\nabla T) = 2\cos\theta
H + \left( (\log f)'\circ\iota \right) (1+\cos^2\theta).
\end{equation}
Remark that this formula yields the following.
See also Lemma 3.1 and Corollary 3.2 in \cite{kn:Ros02}.
\begin{proposition}
There are no compact minimal surfaces in $I\times_f\E^2$ if $
f$ is monotonic.
\end{proposition}
\proof Assume that $(\log f)' \geq 0$ (resp. $\leq 0$) and that $M$
is a compact, minimal surface in $I\times_f\E^2$. By integrating
\eqref{eq:laplacian} and taking into account that $H=0$, one obtains
\begin{equation*}
0= \int_M \Delta h \, dM = \int_M \left( (\log f)'\circ\iota \right)
(1+\cos^2\theta) \, dM \geq 0 \ ({\rm resp.\ } \leq 0).
\end{equation*}
It follows that $(\log f)'\circ\iota = 0$, that is $f$ is constant on
$M$ and the proposition follows immediately.
\endproof

We now consider non-minimal constant angle
surfaces with harmonic height function. If $\sin \theta =0$, then $h$ is constant.
If $\cos\theta=0$, then \eqref{eq:laplacian} implies that $f'=0$ on $M$ such that around $M$
the ambient space is Euclidean and $M$ is a part of a cylinder in the $t$-direction.
If the surface is of type (ii) in Theorem \ref{theo3.1}, with $\theta \in (0,\pi/2)$, then it follows
from \eqref{eq:laplacian} that $f$ is constant on $M$, such that $M$ is part of a plane, hence minimal.
If the surface is of type (i) in Theorem \ref{theo3.1},
with $\theta \in (0,\pi/2)$, then $h$ is harmonic if and only if
\begin{equation}\label{harm}
\sin\theta \cos\theta \lambda + \sigma' = 0.
\end{equation}
From \eqref{harm} and \eqref{lc2} it follows that
\begin{equation*}
\label{harm2}
    \beta=-\cos\theta \frac{1}{\sigma'}.
\end{equation*}
These equations yield that $\beta$ only depends on $u$ and that $\lambda\beta = \frac{1}{\sin\theta}$.
From \eqref{3.1} we easily obtain that $\beta$ is constant.
Therefore $\lambda$ is constant and from \eqref{harm} we obtain that $f(t)=a e^{bt}$.
From \eqref{eq:2} we conclude that the warped product has constant negative sectional curvature.
Without loss of generality we can assume $a=b=1$. It also follows that $\alpha(v)=0$ in \eqref{3.33d} and the surface is given by
\begin{equation*}
\iota(u,v) = \left( u \sin\theta, \cot\theta e^{u\sin \theta} \cos v, \cot\theta e^{u\sin \theta} \sin v \right).
\end{equation*}
Since $\beta$ is constant, $M$ is flat.  Moreover, the
surface is a rotation surface with constant mean curvature
 $H = -(1+\cos^2\theta)/(2\cos\theta)$.
So the surface is a flat constant mean curvature rotation surface in the hyperbolic space.
\begin{remark}
As we have already seen, the ambient $\big(\R^3,\widetilde g=dt^2+e^{2t}(dx^2+dy^2)\big)$ has constant
sectional curvature $-1$. By changing the $t$-coordinate one can obtain the upper half space
model for the hyperbolic 3-space. More precisely, by considering $z=e^{-t}$ one gets
that $(\widetilde M,\widetilde g)$ is isometric to $\big(\H^3_+,g_{-1}\big)$, where
$$
\H^3_+=\left\{(x,y,z)\in\R^3~,~ z>0\right\}\quad,\quad g_{-1}=\frac1{z^2}~\big(dx^2+dy^2+dz^2\big).
$$
In this model, the constant angle surface $M$ obtained above, is given, implicitly by
$(x^2+y^2)z^2=a^2$, $a>0$.
\end{remark}

\end{document}